\documentclass[12pt]{amsart}
\usepackage{pstcol}
\psset{unit=0.7cm}

\advance\textheight1cm 
\def\Ker{{\operatorname{Ker}}}
\def\det{{\operatorname{det}}}
\def\NN{{\Bbb N}}

\def\Hom{{\operatorname{Hom}}}

\def\tame{{\operatorname{tame}}}
\def\Q{\qedsymbol\kern1pt }
\def\SS{\operatorname{S}}
\def\Sur{\operatorname{Sur}}
\def\Aut{{\operatorname{Aut}}}
\def\ggr.aut{{\operatorname{gr.aut}}}

\def\Gr.Hom{{\operatorname{Gr.Hom}}}
\def\L{\operatorname{L}}

\def\Idemp{{\operatorname{Idemp}}}
\def\conv{\operatorname{conv}}
\def\het{\operatorname{ht}}

\def\N{\operatorname{N}}
\def\Spec{\operatorname{Spec}}

\def\gp{\operatorname{gp}}

\def\Ker{{\operatorname{Ker}}}
\def\codim{\operatorname{codim}}

\def\maxSpec{\operatorname{maxSpec}}
\def\Im{\operatorname{Im}}
\def\Pic{\operatorname{Pic}}
\def\Hom{\operatorname{Hom}}
\def\Vect{\operatorname{Vect}}
\def\Pol{\operatorname{Pol}}
\def\GL{\operatorname{GL}}
\def\join{\operatorname{join}}

\def\Col{\operatorname{Col}}
\let\Bbb=\mathbb
\def\RR{{\Bbb R}}

\def\ZZ{{\Bbb Z}}
\def\NN{{\Bbb N}}
\def\TT{{\Bbb T}}

\def\AA{{\Bbb A}}
\def\PP{{\Bbb P}}
\let\cal\mathcal

\let\epsilon=\varepsilon
\let\theta=\vartheta
\let\phi=\varphi

\def\vertex{\circle*{0.15}}

\newtheorem{lemma}{Lemma}[section]
\newtheorem{corollary}[lemma]{Corollary}
\newtheorem{theorem}[lemma]{Theorem}
\newtheorem{proposition}[lemma]{Proposition}

\theoremstyle{definition}

\newtheorem{remark}[lemma]{Remark}

\newtheorem{conjecture}[lemma]{Conjecture}

\input diagrams
\diagramstyle[amstex,PostScript=dvips]
\newarrow{Equal} =====
\newarrow{Into} C--->
\newarrow{Onto} ----{>>}

\begin{document}

\title[Polytopal linear algebra]{Polytopal linear algebra}

\author{Winfried Bruns \and Joseph Gubeladze}
\address{Universit\"at Osnabr\"uck,
FB Mathematik/Informatik, 49069 Osnabr\"uck, Germany}
\email{Winfried.Bruns@mathematik.uni-osnabrueck.de}
\address{A. Razmadze Mathematical Institute, Alexidze St. 1, 380093
Tbilisi, Georgia} \email{gubel@rmi.acnet.ge}
\thanks{The second author was supported by the Max-Planck-Institut f\"ur
 Mathematik at Bonn and TMR grant ERB FMRX CT-97-0107}


\begin{abstract}
We investigate similarities between the category of vector spaces
and that of polytopal algebras, containing the former as a full
subcategory. In Section 2 we introduce the notion of a polytopal Picard group
and show that it is trivial for fields. The coincidence of this group
with the ordinary Picard group for general rings remains an open question.
In Section 3 we survey some of the previous results on the automorphism groups
and retractions. These results support a general conjecture proposed in
Section 4 about the nature of arbitrary homomorphisms of
polytopal algebras. Thereafter a further confirmation of this
conjecture is presented by homomorphisms defined on Veronese singularities.

This is a continuation of the project started in \cite{BG1,BG2,BG3}. The
higher $K$-theoretic aspects of polytopal linear objects will be treated
in \cite{BG4,BG5}.
\end{abstract}

\maketitle

\section{Introduction}
The present work is a continuation of our study of the similarities between
the categories $\Vect(k)$ -- the category of finitely generated vector spaces
over a field $k$, and its natural extension $\Pol(k)$ -- the {\em polytopal
linear category} over $k$, started in the series of papers \cite{BG1,BG2,BG3}.
The category $\Pol(k)$ was first introduced explicitly in \cite{BG3} where we
studied a special class of morphisms, the retractions.

We recall that the objects of $\Pol(k)$ are by definition polytopal
$k$-algebras (discussed in detail in \cite{BGT}), i.~e. the
standard graded $k$-algebras $k[P]$ associated to arbitrary finite convex
lattice polytopes $P\subset\RR^d$ in the following way: the lattice points
$\L_P=\ZZ^d\cap P$ form degree one generators of $k[P]$ and they are subject to
the binomial relations coming from the affine dependencies inside $P$.

Alternatively, $k[P]$ is the semigroup ring $k[S_P]$ of the semigroup
$S_P\subset\ZZ^d$ generated by $\{(x,1)\mid x\in\L_P\}$.

The embedding $\Vect(k)\subset\Pol(k)$, resulting from viewing a vector space
$V$ as the degree one component of its symmetric algebra
$\SS_k(V)=k[X_1,\dots,X_{\dim_kV}]$, makes $\Vect(k)$ a full subcategory of
$\Pol(k)$. Obviously, the latter category is far from being additive but it
reveals many surprising similarities with $\Vect(k)$.

This work provides further results supporting the analogy. In particular, in
Section 2 we show that {\em polytopal Picard groups} defined as the groups of
certain autoequivalences of $\Pol(k)$, are trivial (i.~e.\ coincide with
$\Pic(k)$). We work exclusively over a field $k$, which is often assumed to be
algebraically closed. But in Remark \ref{r2.2} we indicate an approach that
enables one to use arbitrary commutative rings. In the general definition of
the categories $\Vect(k)$ and $\Pol(k)$ the hom-sets are replaced by
appropriate affine schemes. This definition must already be used for fields in
general. The description of $\Pic^{\Pol}(R)$ for general
(commutative) rings $R$ remains an open question.

In order to present a coherent picture of polytopal linear algebra and to ease
references throughout the text, we recall some of the results from \cite{BG1}
and \cite{BG2} in Section 3; they concern the automorphism groups and the
retractions in $\Pol(k)$. In Section 4 we propose a conjecture describing
arbitrary homomorphisms in $\Pol(k)$. Roughly, it says that the homomorphisms
are obtained from the trivial ones by a sequence of standard procedures encoded
in the shapes of the underlying polytopes. In particular, the arithmetic of $k$
is irrelevant in the description of $\Pol(k)$ because everything is determined
on the combinatorial level.

The results obtained so far \cite{BG1,BG2} can be viewed as a confirmation of
refined versions of this conjecture for special classes of morphisms, namely
automorphisms and retractions.

Thereafter in Section 4 we provide further evidence towards our conjecture by
an explicit description of homomorphisms from Veronese subalgebras of
polynomial rings. This result, in conjunction with the results from \cite{BG2},
provides a complete description of the variety of idempotent endomorphisms of
$k[P]$ when $k$ is algebraically closed and $P$ is a lattice polygon (i.~e.\ $\dim P=2$).

The approach developed in \cite{BG3} suggests a further generalization to the
even more general category of {\em polyhedral algebras} and their graded
homomorphisms. This corresponds to the passage from single polytopes to {\em
lattice polyhedral complexes} in the sense of \cite{BG3}. However, in this
article we do not pursue such level of generality and only remark that even the
subclass of simplicial complexes (i.~e.\ the category of Stanley-Reisner rings
and their graded homomorphisms) provides interesting possibilities for the
generalization of linear algebra.

The arguments in Section 2 below use results presented in Section 3. But we
resort to this order of exposition in analogy with the classical hierarchy --
Picard groups, retractions, automorphisms. The objects just listed constitute
the subject of `classical' algebraic $K$-theory (Bass \cite{Ba2}). In
\cite{BG4,BG5} we consider higher $K$-theoretical aspects of the polytopal
generalization of vector spaces.

For the standard terminology in category theory we refer to MacLane \cite{Ma}.

As usual, `standard graded $k$-algebra' means a graded $k$-algebra
$k\oplus A_1\oplus A_2\oplus\cdots$, generated by $A_1$. In what
follows `homomorphism' always means `graded homomorphism'.

For a semigroup $S$ its group of differences (the universal group of $S$)
will be denoted by $\gp(S)$.

Finally, we are grateful to the referee for pointing out to us the references
\cite{Fr} and \cite{Wa}.

\section{Polytopal Picard groups}

Assume $k$ is a field. Then the group of covariant
$k$-linear autoequivalences of $\Vect(k)$, modulo functor isomorphisms,
is a trivial group. Here a functor $F:\Vect(k)\to\Vect(k)$ is called
`$k$-linear' if the mappings
$$
\Hom(V,W)\to\Hom(F(V),F(W)),\quad V,W\in|\Vect(k)|,
$$
are k-linear homomorphisms of vector spaces. This triviality follows from
the fact that
the mentioned group is naturally isomorphic to the Picard group
$\Pic(k)$ (=0) -- an observation valid for any commutative ring $R$. More
precisely, the assignments $F\mapsto F(R)$ and $L\mapsto L\otimes-$ for
$F:{\Bbb M}(R)\to{\Bbb M}(R)$ and $L\in\Pic(R)$ establish an isomorphism
between the group of $R$-linear covariant autoequivalences
of ${\Bbb M}(R)$ -- the category of finitely generated $R$-modules, modulo
functor isomorphisms, and $\Pic(R)$ -- the group of invertible
$R$-modules up to isomorphism \cite{Ba1}.

If $k$ is an algebraically closed field then the condition on a functor
$F:\Vect(k)\to\Vect(k)$ to be $k$-linear is equivalent to the requirements
that the mappings between affine spaces
$$\Hom(V,W)\to\Hom(F(V),F(W)),\quad V,W\in|\Vect(k)|,
$$
are algebraic and, simultaneously, $k^*$-equivariant with respect to the
action
$$
k^*\times\Hom(V,W)\mapsto\Hom(V,W),\quad (a,\phi)\mapsto a\phi.
$$
In the category $\Pol(k)$ both these requirements make sense
(under the assumption $k$ is algebraically closed). In fact, the
sets $\Hom(k[P],k[Q])$ carry natural $k$-variety structures as
follows. An element $f\in\Hom(k[P],k[Q])$ can be identified with
the corresponding matrix
$$
(a_{ij})\in M_{m\times n}(k),\ \ f(x_i)=\sum_{j=1}^na_{ij}y_j,\qquad
x_i\in\L_P,\ y_j\in\L_Q.
$$
Then the equations, defining the Zariski closed subset
$$
\Hom(k[P],k[Q])\subset\AA_k^{mn},\qquad m=\#\L_P,\ n=\#\L_Q,
$$
are derived by the following procedure.
The binomial relations between the $x_i$ are preserved by the $f(x_i)\in
k[Q]_1$ (the degree 1 component of $k[Q]$). After passing to the canonical
$k$-linear expansions as linear forms of monomials in the $y_j$ and comparing
corresponding coefficients (at this point the binomial dependencies between the
$y_j$ are used) we get the desired system of homogeneous equalities $F_s=0$,
$F_s\in k[X_1,\ldots,X_{mn}]$.

\

{\it Observation.} Clearly, we could derive similarly certain homogeneous
polynomials $G_t\in\ZZ[X_1,\ldots,X_{mm}]$ by substituting $\ZZ$ for the
field $k$. Then the polynomials $F_s$ are just specializations of the
$G_t$ under the
canonical ring homomorphism $\ZZ[X_1,\ldots,X_{mn}]\to k[X_1,\ldots,X_{mn}]$.
In particular, the varieties $\Hom(k[P],k[Q])$ are defined over $\ZZ$ and the
corresponding defining integral equations only depend on the polytopes $P$ and
$Q$.

\

As for the $k^*$-equivariant structure, we observe that
any object $A\in|\Pol(k)|$ is naturally equipped with the following
$k^*$-action:
$$
(a_0\oplus a_1\oplus a_2\oplus\cdots)^{\xi}=
a_0\oplus\xi a_1\oplus\xi^2a_2\oplus\cdots,
$$
which induces the algebraic action
$$
k^*\times\Hom(A,B)\to\Hom(A,B)
$$
as follows:
$$
f^{\xi}(a)=f(a^{\xi}),\qquad f\in\Hom(A,B),\ \xi\in k^*,\ a\in A.
$$
Notice that we obtain the same action on $\Hom(A,B)$ by requiring
$$
f^{\xi}(a)=\xi f(a)\ \text{for all}\ f\in\Hom(A,B),\ \xi\in k^*,\
a\in A_1.
$$
It is natural to ask whether the group  of the covariant autoequivalences
(up to functor isomorphism) of  $\Pol(k)$, for which the mappings
$$
\Hom(A,B)\to\Hom(F(A),F(B)),\quad A,B\in|\Pol(k)|
$$
respect the $k$-variety structures and are $k^*$-equivariant, is trivial.
For short, is the `polytopal Picard group' $\Pic^{\Pol}(k)$ trivial?

\begin{remark}\label{r2.2}
We can define the category $\Pol(R)$ for any ring $R$ as
follows. It is the category enriched on the (symmetric) monoidal
category of affine $\Spec(R)$-schemes whose objects are the polytopal
algebras over $R$ and the  hom-schemes $\Hom(R[P],R[Q])$ are the affine
schemes $\Spec(R^{mn}/(G_t))$. Here the $Q_t$ are the same polynomials as
in the observation above. (For the generalities on enriched categories see
\cite{Du}.)

In order to simplify the notation we will consider the underlying rings
instead of the affine schemes. Thus $\Hom(R[P],R[Q])= R^{mn}/(G_t)$. The
equivariant structure on the hom-schemes is encoded into the ring
homomorphisms
$$
\Hom(R[P],R[Q])\to\Hom(R[P],R[Q])[X,X^{-1}],\qquad a\mapsto Xa,
$$
and the composition operation is given by the naturally defined ring
homomorphisms
$$
\Hom(R[P],R[L])\to\Hom(R[P],R[Q])\otimes\Hom(R[Q],R[L]).
$$
Clearly, if $R=k$ is an algebraically closed field, the two definitions of
$\Pol(R)$ are equivalent.

Now we can define $\Pic^{\Pol}(R)$ as the group of those covariant
autoequivalences of $\Pol(R)$ which on the hom-rings induce ring homomorphisms
respecting the $\Spec(\ZZ[X,\allowbreak X^{-1}])$-equivariant structures. In
particular, we have defined $\Pic^{\Pol}(k)$ for general fields.

The proof we present below for algebraically closed fields yields the equality
$\Pic^{\Pol}(k)=0$ for arbitrary fields. We leave this to the interested reader
and only remark that the crucial fact is that the main result of \cite{BG1}
(Theorem \ref{t3.2} below) has been proved for general fields.

The lack of an analogous description of the automorphism groups over a general
ring of coefficients is the obstacle in describing the group $\Pic^{\Pol}(R)$.
We expect that this is a trivial group, being a polytopal counterpart of the
group of $R$-linear autoequivalences (up to functor isomorphism) of the
category of finitely generated free $R$-modules -- a trivial group.
\end{remark}

\begin{remark}\label{r2.3}
The last step in the proof of Theorem \ref{t2.4} uses the fact from \cite{BGT}
that for a polytope $P$ and a field $k$ the polytopal algebra $k[cP]$ is
quadratically defined (i.~e.\ by degree 2 equations) whenever $c\geq\dim P$.
This however does not create an additional difficulty in generalizing the
result on polytopal Picard groups from fields to arbitrary rings. It is
an elementary fact that the condition on a polytopal ring to be
quadratically defined depends only on the combinatorial structure of
the polytope.
\end{remark}

In the remaining part of this section $k$ is an algebraically closed field.

\begin{theorem}\label{t2.4}
$\Pic^{\Pol}(k)=0$.
\end{theorem}

\begin{proof} Let $[F]\in\Pic^{\Pol}(k)$. We want to show that there are
isomorphisms $\bigl(\sigma_A:A\to F(A)\bigr)_{|\Pol(k)|}$ such that
$$
F(f)=\sigma_B\circ f\circ\sigma_A^{-1}
$$
for every homomorphism $f:A\to B$ in $\Pol(k)$.

First of all notice that we can work on an arbitrarily fixed skeleton of
$\Pol(k)$: all the notions we are dealing with are invariant under such a
passage. Henceforth $\Pol(k)$ is the fixed skeleton. By \cite{Gu} one knows
that $k[P]$ determines (up to an affine integral isomorphism) the polytope $P$
(this is so even in the category of all commutative $k$-algebras). Therefore,
we assume that for each object $A\in|\Pol(k)|$ there is a unique polytope $P$
such that $A=k[P]$ and different polytopal algebras determine non-isomorphic
polytopes. By a suitable choice of the skeleton we can also assume that the
objects of $\Vect(k)$ are of the type $k[\Delta_n]$ for the unit $n$-simplices
$\Delta_n$, $n=-1,0,1,2,\dots$ (by convention $k[\Delta_{-1}]=k$\pagebreak[3]). Also, we
will use the notation $k[t]=k[\Delta_0]$.
\bigskip

\noindent{\em Step 1.} We claim that
\begin{itemize}
\item[(i)]
$\dim A=\dim F(A)$ (Krull dimension),
\item[(ii)]
$\dim_kA_1=\dim_kF(A)_1$ ($k$-ranks of the degree 1 components),
\item[(iii)]
$F(f)$ is surjective (the degree 1 component $F(f)_1$ is injective) if and only if $f$ is surjective ($f_1$ is injective).
\end{itemize}

In fact, it follows from Theorem \ref{t3.2} below that for any polytopal algebra
$k[P]$ its Krull dimension is the dimension of a maximal torus of
the linear subgroup
$$
\Gamma_k(P)=\Aut_{\Pol(k)}(k[P])\subset\Aut_k{k[P]_1};
$$
it is certainly an invariant of $F$ -- hence (i).

The claims (ii) and (iii) follow from the observation that
both the injectivity and surjectivity conditions can be reformulated
in purely categorical terms by using morphisms originating from $k[t]$.
\bigskip

\noindent {\em Step 2}. Observe that $F$ restricts to an autoequivalence of
$\Vect(k)$. This follows from the claims (i) and (ii) in Step 1 and the fact
that polynomial algebras are the only polytopal algebras whose Krull dimensions
coincide with the $k$-ranks of the degree 1 components.

Next we correct $F$ on $\Vect(k)$ in such a way that $F|_{\Vect(k)}= {\bf
1}_{\Vect(k)}$.

We know that $F|_{\Vect(k)}\approx{\bf 1}_{\Vect(k)}$. This means that there
are elements $\rho_V\in\Aut(V)$ for all $V\in|\Vect(k)|$ such that
$$
\forall\ U,V\in|\Vect(k)|\ \ \ (f:U\to V)\mapsto(\rho_V\circ
f\circ\rho_U^{-1}:U\to V).
$$
For each $A\in|\Pol(k)|\setminus|\Vect(k)|$ put $\rho_A={\bf 1}_A$ and let
$H:\Pol(k)\to\Pol(k)$ be the autoequivalence determined as follows: it is the
identity on the objects and
$$
H(f)=\rho_B^{-1}\circ f\circ\rho_A
$$
for a morphism $f:A\to B$ in $\Pol(k)$. Then $H$ represents the neutral element
of $\Pic^{\Pol}(k)$. Now the functor $G=H\circ F$ is isomorphic to $F$ and it
restricts to the identity functor on $\Vect(k)$.

Without loss of generality we can therefore assume $F|_{\Vect_k}={\bf
1}_{\Vect_k}$. \bigskip

\noindent {\em Step 3.} For any positive dimensional object $k[P]\in|\Pol(k)|$
we fix a bijective mapping $\Delta_P$ from the set of vertices of
$\Delta_{n-1}$, $n=\#\L_P$ to the set of lattice points of $P$. The resulting
surjective homomorphism $k[\Delta_{n-1}]\to k[P]$ will also be denoted by
$\Delta_P$.

Every matrix $\alpha\in\GL_n(k)$ gives rise to a graded $k$-surjective
homomorphism
$$
\alpha\Delta_P:k[\Delta_{n-1}]\to k[P]
$$
whose degree 1 component is given by
$$
(\alpha\Delta_P)_1=(\Delta_P)_1\circ\alpha_*,
$$
where $\alpha_*$ is the linear transformation of the $k$-vector space
$k[\Delta_{n-1}]_1$ determined by $\alpha$, and $(\Delta_P)_1$ is the degree 1
component of $\Delta_P$.

By Step 2 $F(k[\Delta_{n-1}])=k[\Delta_{n-1}]$ for all $n\in\NN$. Therefore,
for a polytope $P$, satisfying the condition $F(k[P])=k[P]$, the claim (iii) in
Step 1 implies the following commutative square
$$
\begin{diagram}
\Hom(k[\Delta_{n-1}],k[P])&\rTo^F&\Hom(k[\Delta_{n-1}],k[P])\\
\uInto                    &      &\uInto                    \\
\Sur(k[\Delta_{n-1}],k[P])&\rTo^F&\Sur(k[\Delta_{n-1}],k[P])
\end{diagram}
$$
where $\Sur(k[\Delta_{n-1}],k[P])$ denotes the set of surjective homomorphisms.
By the definition of $\Pic^{\Pol}(k)$ the horizontal mappings are
$k^*$-equivariant automorphisms of $k$-varieties.

We have the following obvious equalities:
\begin{align*}
\Hom(k[\Delta_{n-1}],k[P])&=\{\gamma\Delta_P|\ \gamma\in M_{n\times n}(k)\},\\
\Sur(k[\Delta_{n-1}],k[P])&=\{\gamma\Delta_P|\ \gamma\in\GL_n(k)\}.
\end{align*}
After the appropriate identifications with $M_{n\times n}(k)$ and $\GL_n(k)$
respectively we arrive at the commutative square of $k$-varieties
$$
\begin{diagram}
M_{n\times n}(k) & \rTo^F & M_{n\times n}(k)\\
\uInto           &        &\uInto\\
\GL_n(k)         &\rTo^F  &\GL_n(k)
\end{diagram}
$$
whose horizontal arrows are algebraic $k^*$-equivariant automorphisms for the
diagonal $k^*$-actions. Thus the upper horizontal mapping is a linear
non-degenerate transformation of the $k$-vector space $M_{n\times n}(k)$,
which leaves the subset $\GL_n(k)\subset M_{n\times n}(k)$
invariant, i.~e. the matrix degeneracy locus $M_{n\times
n}(k)\setminus\GL_n(k)$ is invariant under $F$. Then by Proposition \ref{c3.3}
below there are only two possibilities:
\begin{itemize}
\item[(a)]
there are $\alpha,\beta\in\GL_n(k)$ such that $F(\gamma\Delta_P)
=(\beta\gamma\alpha)\Delta_P$ for all $\gamma\in M_{n\times n}(k)$, or
\item[(b)]
there are $\alpha,\beta\in\GL_n(k)$ such that
$F(\gamma\Delta_P)=(\beta\gamma^{T}\alpha)\Delta_P$ for all $\gamma\in
M_{n\times n}(k)$, where $-^T$ is the transposition.
\end{itemize}
Consider the commutative square
$$
\begin{diagram}
k[\Delta_{n-1}]       &\rTo^{\gamma_*}&k[\Delta_{n-1}]\\
\dTo^{\gamma\Delta_P} &               &\dTo_{\Delta_P}\\
k[P]                  &\rEqual        &k[P],
\end{diagram},
$$
where $\gamma\in M_{n\times n}(k)$ is arbitrary matrix and the degree 1
component of the upper horizontal mapping is the linear transformation given by
$\gamma$. By applying the functor $F$ to this square and using the equality
$F|_{\Vect(k)}={\bf 1}_{\Vect(k)}$ (Step 2) we arrive at the following
equalities in the corresponding cases:
$$
\text{(a)}\quad
(\beta\gamma\alpha)\Delta_P=(\beta\alpha)\Delta_P\circ\gamma_*,\qquad
\text{(b)}\quad
(\beta\gamma^T\alpha)\Delta_P=(\beta\alpha)\Delta_P\circ\gamma_*.
$$
Identifying $\L_{\Delta_{n-1}}$ and $\L_P$ along $\Delta_P$ we get the matrix
equalities:
$$
\text{(a)}\quad \beta\gamma\alpha=\gamma\beta\alpha,\qquad
\text{(b)}\quad \beta\gamma^T\alpha=\gamma\beta\alpha.
$$
First notice that case (b) is excluded, i.~e.\ there is no matrix
$\beta\in\GL_n(k)$ for which the following holds:
$$
\forall\ \gamma\in M_{n\times n}(k)\quad\beta\gamma\beta^{-1}=\gamma^T.
$$
This follows from running  $\gamma$ through the set of standard basic
matrices (i.~e. the matrices with only one entry 1 and 0s elsewhere).

For case (a) we have
$$
\forall\ \gamma\in M_{n\times n}(k)\quad \beta\gamma\beta^{-1}=\gamma.
$$
Then $\beta$ is in the center of $M_{n\times n}(k)$ (in particular, it is a
scalar matrix). So we can write
$$
F(\gamma\Delta_P)=(\gamma\beta\alpha)\Delta_P
$$
We arrive at the
\medskip

\noindent {\bf Claim.} For each lattice polytope $P$ with $F(k[P])=k[P]$, there
is $\alpha_P\in\GL_n(k)$, $n=\#\L_P$, such that
$$
\forall\ \gamma\in\ M_{n\times n}(k)\quad
F(\gamma\Delta_P)=(\gamma\alpha_P)\Delta_P.
$$
\bigskip

\noindent {\em Step 4.} Now we show that for every polytope $P$ the linear
automorphism $(\alpha_P)_*$ of the $k$-vector space $k[P]_1$, determined by the
matrix $\alpha_P$ (Step 3), belongs to the closed subgroup
$$
\Gamma_k(P)\subset\Aut_k(k[P]_1),
$$
provided $F(k[P])=k[P]$.

The functor $F$ induces a $k^*$-equivariant automorphism of the $k$-variety
$\Hom(k[P],\allowbreak k[t])$. On the other hand we have the natural
identification
$$
\Hom(k[P],k[t])=\maxSpec(k[P]),
$$
where the right hand side denotes the variety of closed points of
$\Spec(k[P])$. Therefore, there exists an automorphism $\psi$ of the
$k$-algebra $k[P]$ -- {\em not a priori} graded, such that the mapping
$$
F:\Hom(k[P],k[t])\to\Hom(k[P],k[t])
$$
is given by $\phi\mapsto\phi\circ\psi$, and, moreover, it is $k^*$-equivariant.
It follows that $\psi$ is a $k^*$-equivariant automorphism of $k[P]$. It is
easily seen that a $k^*$-equivariant automorphism is graded (and conversely).
Therefore $\psi_1\in\Gamma_k(P)$. (Here we identify elements of $\Gamma_k(P)$
with their degree one components, which are linear automorphisms of
$\bigoplus_{\L_P}k$.)

We let $P^*:k[P]\to k[t]$ denote the homomorphism which sends $\L_P$ to $t$.
For any toric automorphism $\tau\in\TT_k(P)$ (i.e. an automorphism for which
any element of $\L_P$ is an eigenvector, $\TT_k(P)$ is the group of such
automorphisms, see Section 3) we have the commutative diagram
$$
\begin{CD}
k[\Delta_{n-1}]@>{(\Delta_{n-1})^*}>>k[t]\\
@V{\tau^*\Delta_P}VV@AA{P^*}A\\
k[P]@>>{\tau^{-1}}>k[P]
\end{CD}
$$
where $\tau^*$ refers to the diagonal $n\times n$-matrix corresponding to
the degree 1 component of $\tau$.

In view of what has been said above and of Step 3, an application of $F$ to the
last commutative diagram yields the equality
$$
(\Delta_{n-1})^*=(P^*\circ\tau^{-1}\circ\psi)\circ
{\big (}(\tau^*\alpha_P)\Delta_P{\big )}.
$$
Let $\psi^*\in M_{n\times n}(k)$ be the matrix of the degree 1 component of
$\psi$, and put $\omega=\alpha_P\psi^*$. Then the equality can be reformulated
into the condition:
\begin{itemize}
\item[($*$)]
for any $\tau\in\TT_k(P)$ the sum of the entries of each row in the matrix
$\tau^*\omega(\tau^*)^{-1}$ is 1, that is
$$
\bigl((\tau^{-1}\circ\omega_*\circ\tau)(x_i)= {\textstyle \sum_{j=1}^n}
c_{ij}x_j\bigr) \implies\bigl({\textstyle \sum_{j=1}^n} c_{ij}=1\big ),
$$
where $\omega_*$ is the linear transformation of $k[P]_1$ determined by
$\omega$, and $\L_P=\{x_1,\dots,x_n\}$.
\end{itemize}

Now we derive from ($*$) that $\omega_*={\bf 1}_{k[P]_1}$, that is
$(\alpha_P)_*=\psi^{-1}\in\Gamma_k(P)$. Fix $i\in[1,n]$ and $\tau\in\TT_k(P)$.
We put
$$
\omega_*(x_i)=\sum_{j=1}^na_jx_j
$$
and
$$
\tau(x_j)=b_jx_j,\qquad j\in[1,n].
$$
Then
$$
(\tau^{-1}\circ\omega_*\circ\tau)(x_i)=\sum_{j=1}^n\frac{b_ia_j}{b_j}x_j.
$$
Consider the Laurent polynomial
$$
l_i=\sum_{j=1}^na_j\frac{x_j}{x_i}\in k[\ZZ^d]\ (
=k[z_1,z_1^{-1},\dots,z_d,z_d^{-1}]).
$$
Without loss of generality we assume that $d=\dim P$ and $\gp(S_P)=\ZZ^{d+1}$.

Observe that the assignment
$$
\frac{x_j}{x_i}\mapsto\frac{b_ia_j}{b_j}\cdot\frac{x_j}{x_i},\qquad j\in[1,n],
$$
extends uniquely to a toric automorphism of $k[\ZZ^d]$ (the torus $(k^*)^d$
acts tautologically on its coordinate ring $k[\ZZ^d]$). Conversely, any toric
automorphism of $k[\ZZ^d]$ can be obtained in this way from some element of
$\TT_k(P)$. Therefore
$$
\bigl\{\xi(l_i)\mid \xi\in(k^*)^d\bigr\}=
\bigl\{x_i^{-1}((\tau^{-1}\circ\omega_*\circ\tau)(x_i))\mid
\tau\in\TT_k(P)\bigr\}.
$$
In particular, the sum of the coefficients of the Laurent polynomial $\xi(l_i)$
is 1 for any $\xi\in(k^*)^d$, in other words
$$
l_i(\xi_1,\dots,\xi_d)=1
$$
for all $\xi_1,\dots,\xi_d\in k^*$. Taking into account the infinity of
$k$ we conclude $l_i=1$, that is $\omega_*={\bf 1}_{k[P]_1}$.
\bigskip

\noindent {\em Step 5.} Let $\Pol(k)'\subset\Pol(k)$ denote the full
subcategory whose objects are those polytopal algebras $k[P]$ for which
$F(k[P])=k[P]$. Now we show that there is $G\in\Pic^{\Pol}(k)$ such that
$F\approx G$ and $G|_{\Pol(k)'}={\bf 1}_{\Pol(k)'}$.

As in Step 2, this means that we have to show the existence of elements
$$
\sigma_P\in\Gamma_k(P),\qquad k[P]\in|\Pol(k)'|,
$$
such that $F(f)_1=\sigma_Q\circ f_1\circ\sigma_P^{-1}$ for any morphism
$f:k[P]\to k[Q]$ in $\Pol(k)'$. Consider the commutative diagram
$$
\begin{CD}
k[\Delta_{m-1}]@>{\tilde f}>>k[\Delta_{n-1}]\\
@V{\Delta_P}VV@VV{\Delta_Q}V\\
k[P]@>>f>k[Q]
\end{CD}
$$
where $f$ is any morphism in $\Pol(k)'$, $m=\#\L_P$, $n=\#\L_Q$,
and $\tilde f$ is the unique lifting of $f$. By Step 3 we
know that $F$ transforms this square into the square
$$
\begin{CD}
k[\Delta_{n-1}]@>{\tilde f}>>k[\Delta_{m-1}]\\
@V{\alpha_P\Delta_P}VV@VV{\alpha_Q\Delta_Q}V\\
k[P]@>>{F(f)}>k[Q]
\end{CD}
$$
Therefore, looking at the degree 1 components we conclude
$$
F(f)_1=(\alpha_Q)_*\circ f_1\circ(\alpha_P)_*^{-1}.
$$
So by Step (4) the system
$$
\bigl\{(\alpha_P)_*\mid k[P]\in|\Pol(k)'|\bigr\}
$$
is the desired one.
\bigskip

\noindent {\em Step 6.} By the previous step we can assume that
$F|_{\Pol(k)'}={\bf 1}_{\Pol(k)'}$. Now we complete the proof by showing that
this assumption implies $\Pol(k)'=\Pol(k)$, and therefore $F={\bf
1}_{\Pol(k)}$.

Assume $P$ is a lattice polytope. We let $\Q$ denote the unit square and
consider all the possible integral affine mappings:
$$
{\cal M}=\{\mu:\Q\to P\}.
$$
We define the diagram ${\Bbb D}(k,P)$ as follows. It consists of
\begin{itemize}
\item[(1)]
$\#(\L_P)$ copies of $k[t]$ indexed by the elements of
$\L_P=\{x_1,\dots,x_n\}$,
\item[(2)]
$\#({\cal M})$ copies of $k[\Q]$: $\{k[\Q]_\mu\}_{\cal M}$,
\end{itemize}
and the following morphisms between them
$$
k[t]_{x_i}\to k[\Q]_\mu, \qquad t\mapsto z\in\L_{\Q},\ \ \text{if}\ \
\mu(z)=x_i.
$$
We claim that
$$
\varinjlim{\Bbb D}(k,P)=k[P]
$$
whenever the defining ideal of the toric ring $k[P]$ is generated by
quadratic binomials, where the direct limit is considered in the category
of all (commutative) $k$-algebras.

In fact, it is clear that the mentioned limit is always a standard graded
$k$-algebra, whose degree 1 component has $k$-dimension equal to $\#\L_P$, and
that there is a surjective graded $k$-homomorphism
$$
{\Bbb L}(k,P)=\varinjlim{\Bbb D}(k,P)\to k[P].
$$
This is so because of the cone over the diagram ${\Bbb D}(k,P)$ with
vertex $k[P]$ where $t\in k[t]_{x_i}$ is mapped to $x_i$ and the
vertices of $\Q$ from $k[\Q]_{\mu}$ are mapped accordingly to $\mu$.
Therefore, we only need to make sure that $t_it_j-t_kt_l=0$ in ${\Bbb
L}(k,P)$ whenever $x_ix_j-x_kx_l=0$ in $k[P]$, where $t_i$ is the image
of $t\in k[t]_{x_i}$ in ${\Bbb L}(k,P)$, and similarly for $t_j$, $t_k$
and $t_l$.  But if $x_ix_j-x_kx_l=0$ in $k[P]$ then these four points
in $P$ belong to $\Im(\mu)$ for some $\mu\in{\cal M}$ and the desired
equality is encoded into the diagram ${\Bbb D}(k,P)$.

Having established the equality above for quadratically defined polytopal rings
we now show that $k[cP]\in|\Pol(k)'|$ for every lattice polytope $P$ and all
$c\geq\dim P$. (Here $cP$ denotes the $c$th homothetic multiple of $P$.)

By Theorem 1.3.3(a) in \cite{BGT} the defining ideal of the toric ring $k[cP]$
is generated by quadratic binomials for $c\geq\dim P$. Therefore,
$$
{\Bbb L}(k,cP)=k[cP]
$$
for such $c\in\NN$. Now observe that the claims (i) and (ii) in Step 1
imply that either $F(k[\Q])=k[\Q]$ or $F(k[\Q])=k[T]$ for a lattice triangle
$T$ with 4 lattice points.

Let us show that $F(k[\Q])$ cannot be a `triangle ring'. Up to isomorphism
there are only two lattice triangles with 4 lattice points
$$
\conv{\big(}(1,0),(0,1),(-1,0){\big)}\qquad\text{and}\qquad
\conv{\big(}(1,0),(0,1),(-1,-1){\big)}
$$
\begin{figure}[htb]
\begin{center}
\psset{unit=1cm}
\begin{pspicture}(-0.4,-1)(10.5,1)
\def\vertex{\pscircle[fillstyle=solid,fillcolor=black]{0.07}}
\psset{linewidth=0pt}
\pspolygon[fillstyle=solid,fillcolor=yellow](0,0)(0,1)(1,1)(1,0)
\pspolygon[fillstyle=solid,fillcolor=yellow](4,0)(6,0)(5,1)
\pspolygon[fillstyle=solid,fillcolor=yellow](8,-1)(10,0)(9,1)
\psset{linewidth=0.7pt,linecolor=black}
\psline{->}(5,1)(4.05,0.05)
\psline{->}(5,1)(5,0.05)
\psline{->}(5,1)(5.95,0.05)
\psline{->}(5,0)(4.05,0)
\psline{->}(5,0)(5.95,0)
\psline{->}(0,0)(0.95,0)
\psline{->}(0,0)(0,0.95)
\psline{->}(1,1)(1,0.05)
\psline{->}(1,1)(0.05,1)
\pspolygon(8,-1)(10,0)(9,1)
\multirput(0,0)(1,0){2}{\vertex}
\multirput(0,1)(1,0){2}{\vertex}
\multirput(4,0)(1,0){3}{\vertex}
\rput(5,1){\vertex}
\rput(8,-1){\vertex}
\rput(9,1){\vertex}
\rput(10,0){\vertex}
\rput(9,0){\vertex}
\end{pspicture}
\end{center}
\caption{}\label{Fig1}
\end{figure}
By Theorem \ref{t3.2} below the triangle $T$ with the property $F(k[\Q])=k[T]$
must have the same number of {\em column vectors} as $\Q$ (see Section 3 for
the definition). But this is not the case because the second triangle
has no column vectors whereas the first has 5 of them, and the square only 4.
(The column vectors are indicated in the figure.) Therefore
$k[\Q]\in|\Pol(k)'|$ as well.

By the assumption $F|_{\Pol(k)'}={\bf 1}_{\Pol(k)'}$ we get the natural
homomorphism $k[cP]\to F(k[cP])$ ($c$ as above). Arguing similarly for the
functor $F^{-1}$, we get the natural homomorphism $k[cP]\to F^{-1}(k[cP])$, and
applying $F$ to the latter homomorphism we arrive at the natural homomorphism
$F(k[cP])\to k[cP]$. But every natural (that is compatible with ${\Bbb
D}(k,cP)$) endomorphism $k[cP]\to k[cP]$ must be the identity mapping for
reasons of universality. In particular, $k[cP]$ is a $k$-retract of $F(k[cP])$.
By reasons of dimensions (and claim (i) in Step 1) we finally
get the equality $F(k[cP])=k[cP]$, as required.

Now we are ready to show the equality $\Pol(k)'=\Pol(k)$, completing the
proof of Theorem \ref{r2.3}.

Let $k[P]\in|\Pol(k)|$. Consider two coprime natural numbers
$c,c'\geq\dim P$. We have the commutative square, consisting of
embeddings in $\Pol(k)$,
\begin{equation}\tag{$**$}
\begin{CD}
k[P]@>>>k[cP]\\
@VVV@VVV\\
k[c'P]@>>>k[cc'P]
\end{CD}
\end{equation}
where the horizontal (vertical) mappings send lattice points in
the corresponding polytopes to their homothetic images (centered
at the origin) with factor $c$ and $c'$ respectively. The key
observation is that by restricting to the degree one components we
get the pull back diagram of $k$-vector spaces
$$
\begin{CD}
k[P]_1@>>>k[cP]_1\\
@VVV@VVV\\
k[c'P]_1@>>>k[cc'P]_1.
\end{CD}
$$
This follows from the fact that the following is a pull back diagram of finite sets
$$
\begin{CD}
\L_P@>>>\L_{cP}\\ @VVV@VVV\\ \L_{c'P}@>>>\L_{cc'P}.
\end{CD}
$$

\noindent {\em Caution.} $(**)$ is in general not a pull back diagram of
$k$-algebras. Otherwise, by Theorem 1.3.3 in \cite{BGT}, any polytopal algebra
would be normal, which is not the case. \medskip

That the latter square of finite sets is in fact a pull back diagram
becomes transparent after thinking of it as the (isomorphic) diagram
consisting of the inclusions:
$$
\begin{diagram}
\L_{P(\ZZ^d)}           &\rInto &\L_{P(\frac{\ZZ^d}{c})}\\
\dInto                  &       &\dInto\\
\L_{P(\frac{\ZZ^d}{c'})}&\rInto &\L_{P(\frac{\ZZ^d}{cc'})}
\end{diagram}
$$
where $P({\cal H})$ refers to the lattice polytope $P\subset\RR^d$ with respect
to the intermediate lattice $\ZZ^d\subset{\cal H}\subset\RR^d$. Here again we
assume that $d=\dim P$ and $\gp(S_P)=\ZZ^{d+1}$.

Since $k[cP],k[c'P],k[cc'P]\in\Pol(k)'$ and $F$ is the identity functor on
$\Pol(k)'$, an application of $F$ to the square $(**)$ yields the square of
graded homomorphisms
$$
\begin{CD}
F(k[P])@>>>k[cP]\\ @VVV@VVV\\ k[c'P]@>>>k[cc'P],
\end{CD}
$$
with the same arrows into $k[cc'P]$ as $(**)$. The same arguments as in the
proof of (iii), Step 1 show that the degree 1 component of this square is a
pull-back diagram, in other words $F(k[P])_1=k[P]_1$. Therefore, $F(k[P])_1$
and $k[P]_1$ generate the same algebras, i.~e. $k[P]=F(k[P])$.
\end{proof}

\section{Automorphisms and retractions}
In this section we survey those results of \cite{BG1} and \cite{BG2} that are
related to polytopal linear algebra.

As far as automorphisms are concerned, $k$ will be assumed to be a general, not
necessarily algebraically closed field.

As remarked in Section 2, for a lattice polytope $P\subset\RR^d$ the group
$\Gamma_k(P)=\ggr.aut(k[P])$ is a linear $k$-group in a natural way. It
coincides with $\GL_n(k)$ in the case of the unit $(n-1)$-simplex
$P=\Delta_{n-1}$. The groups $\Gamma_k(P)$ have been named {\em polytopal
linear groups} in \cite{BG1}.

An element $v\in \ZZ^d$, $v\neq 0$, is a {\em column vector} (for $P$) if
there is a facet $F\subset P$ such that $x+v\in P$ for every
lattice point $x\in P\setminus F$. The facet $F$ is called {\em the base
facet} of $v$.
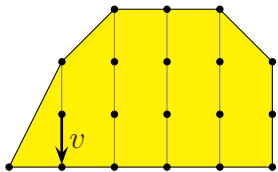
\begin{figure}[htb]
\begin{center}
\begin{pspicture}(0,0)(5,3)
\psset{linewidth=0.4pt}
\def\vertex{\pscircle[fillstyle=solid,fillcolor=black]{0.07}}
\pspolygon[fillstyle=solid,fillcolor=yellow](0,0)(5,0)(5,2)(4,3)(2,3)(1,2)
\psline[linecolor=gray](2,0)(2,3)
\psline[linecolor=gray](3,0)(3,3)
\psline[linecolor=gray](4,0)(4,3)
\psline[linecolor=gray](1,0)(1,2)
\multirput(0,0)(1,0){6}{\vertex}
\multirput(1,1)(1,0){5}{\vertex}
\multirput(1,2)(1,0){5}{\vertex}
\multirput(2,3)(1,0){3}{\vertex}
\rput(1.3,0.5){$v$}
\psline[linewidth=1pt]{->}(1,1)(1,0.05)
\end{pspicture}
\end{center}
\caption{A column structure}\label{Fig2}
\end{figure}
The set of column vectors of $P$ is denoted by $\Col(P)$. A pair $(P,v)$,
$v\in\Col(P)$, is called a {\em column structure}. Let $(P,v)$ be a column
structure and $P_v\subset P$ be the base facet for $v\in\Col(P)$. Then for each
element $x\in S_P$ we set $\het_v(x)=m$ where $m$ is the largest non-negative
integer for which $x+mv\in S_P$. Thus $\het_v(x)$ is the `height' of $x$ above
the facet of the cone $C(S_P)$ corresponding to $P_v$ in direction $-v$. It is
an easy observation that $x+\het_v(x)\cdot v\in S_{P_v}\subset S_P$ for any
$x\in S_P$.

Column vectors are just the dual objects to the roots of the normal fan
$\mathcal N(P)$ in the sense of Demazure (see \cite[Section 3.4]{Oda}).
(For the relationship of Theorem 3.2 with the results of Demazure
\cite{De} and Cox \cite{Cox} on the automorphism groups of complete
toric varieties see \cite[Section 5]{BG1}.)

Let $(P,v)$ be a column structure and $\lambda\in k$. We identify the
vector $v$, representing the difference of two lattice points in $P$,
with the degree $0$ element $(v,0)\in\gp(S_P)\subset k[\gp(S_P)]$. Then
the assignment
$$
x\mapsto (1+\lambda v)^{\het_v(x)}x.
$$
gives rise to a graded $k$-algebra automorphism $e_v^\lambda$ of $k[P]$.
Observe that $e_v^\lambda$ becomes an elementary matrix in the special
case when $P=\Delta_{n-1}$, after the identifications
$k[{\Delta_{n-1}}]=k[X_1,\dots,X_n]$ and $\Gamma_k(P)=\GL_n(k)$. Accordingly
$e_v^\lambda$ is called an {\em elementary automorphism}.

The following alternative description of elementary automorphisms is essential
in showing that all automorphisms in $\Pol(k)$ are {\em tame} (see Remark
\ref{r4.2}(c) below). We may assume $\dim P=d$ and $\gp(S_P)=\ZZ^{d+1}$. By a
suitable integral unimodular change of coordinates we may further assume that
$v=(0,-1,0,\dots,0)$ and that $P_v$ lies in the subspace $\RR^{d-1}$ (thus $P$
is in the upper halfspace). Consider the standard unimodular simplex
$\Delta_{d-1}$ (i.~e. the one with vertices at the origin and the standard
coordinate unit vectors). Clearly, $P$ is contained in a parallel integral
shift of $c\Delta_{d-1}$ for a sufficiently large natural number $c$. Then we
have a graded $k$-embedding $k[P]\to k[c\Delta_{d-1}]$, the latter ring being
just the $c$th Veronese subring of the polynomial ring $k[X_1,\dots,X_d]$.
Moreover, $v=X_1/X_2$. Now the automorphism of $k[X_1,\dots,X_d]$ mapping
$X_2$ to $X_2+\lambda X_1$ and leaving all the other variables fixed induces an
automorphism of $k[c\Delta_{d-1}]$ and the latter restricts to an automorphism
of $k[P]$, which is nothing but the elementary automorphism $e_v^{\lambda}$
above.

As usual, $\AA_k^s$ denotes the additive group of the $s$-dimensional affine
space.

\begin{proposition}\label{p3.1}
Let $v_1,\dots,v_s$ be pairwise different column vectors for $P$ with
the same base facet $F=P_{v_i}$, $i=1,\dots,s$. The mapping
$$
\phi:\AA_k^s\to\Gamma_k(P),\qquad
(\lambda_1,\dots,\lambda_s)\mapsto
e_{v_1}^{\lambda_1}\circ\cdots\circ e_{v_s}^{\lambda_s},
$$
is an embedding of algebraic groups. In particular,
$e_{v_i}^{\lambda_i}$ and $e_{v_j}^{\lambda_j}$ commute for all
$i,j\in\{1,\dots,s\}$ and $(e_{v_1}^{\lambda_1}\circ\cdots\circ
e_{v_s}^{\lambda_s})^{-1} =e_{v_1}^{-\lambda_1}\circ\cdots\circ
e_{v_s}^{-\lambda_s}.$
\end{proposition}

The image of the embedding $\phi$ given by Lemma \ref{p3.1} is denoted by
$\AA(F)$. Of course, $\AA(F)$ may consist only of the identity map of
$k[P]$, namely if there is no column vector with base facet $F$.

Put $n=\dim(P)+1$. The $n$-torus $\TT_n=(k^*)^n$ acts naturally on
$k[P]$ by restriction of its action on $k[\gp(S_P)]$ that is given by
$$
(\xi_1,\dots,\xi_n)(e_i)=\xi_ie_i,\quad
i\in[1,n].
$$
Here $e_i$ is the $i$-th element of a fixed basis of
$\gp(S_P)=\ZZ^n$.  This gives rise to an algebraic embedding
$\TT_n\subset\Gamma_k(P)$, whose image we denote by $\TT_k(P)$.
It consists precisely of those automorphisms of $k[P]$ which
multiply each monomial by a scalar from $k^*$.

The (finite) automorphism group $\Sigma(P)$ of the semigroup $S_P$ is
also a subgroup of $\Gamma_k(P)$. It is exactly the group of automorphisms
of $P$ as a lattice polytope.

Next we recapitulate the main result of \cite{BG1}. It should be viewed a
polytopal generalization of the standard linear algebra fact that any
invertible matrix over a field can be reduced to a diagonal matrix (generalized
to toric automorphisms in the polytopal setting) using elementary
transformations on columns (rows). Moreover, we have normal forms for such
reductions reflecting the fact that the elementary transformations can be
carried out in an increasing order of the column indices.

\begin{theorem}\label{t3.2}
Let $P$ be a convex lattice $n$-polytope and $k$ a field.
Every element $\gamma\in\Gamma_k(P)$ has a (not uniquely determined)
presentation
$$
\gamma=\alpha_1\circ\alpha_2\circ\cdots\circ\alpha_r\circ\tau\circ\sigma,
$$
where $\sigma\in\Sigma(P)$, $\tau\in\TT_k(P)$, and
$\alpha_i\in\AA(F_{i})$ such that the facets $F_i$ are pairwise
different and $\#\L_{F_i}\le \#\L_{F_{i+1}}$, $i\in[1,r-1]$.

We have $\dim\Gamma_k(P)=\#\Col(P)+n+1$ (the left hand side is the Krull
dimension of the group scheme $\Gamma_k(P)$), and $\TT_k(P)$ is a maximal torus
in $\Gamma_k(P)$, provided $k$ is infinite.
\end{theorem}

As an application beyond the theory of toric varieties we mention
that Theorem \ref{t3.2} provides yet another proof of the
classical description of the graded automorphisms of determinantal
rings -- a result which goes back to Frobenius \cite[p.~99]{Fr}
and has been re-proved many times since then. See, for instance,
\cite{Wa} for a group-scheme theoretical approach which involves
general commutative rings of coefficients and the classes of
generic symmetric and alternating matrices.

\begin{proposition}\label{c3.3}
Let $k$ be any field, $X$ an $m\times n$ matrix of indeterminates, and
$R=k[X]/I_{r+1}(X)$ the residue class ring of the polynomial ring $k[X]$
in the entries of $X$ modulo the ideal generated by the $(r+1)$-minors
of $X$, $1\leq r<\min(m,n)$. Let $G=\ggr.aut(R)$ and $G^0$ denote the image of
the mapping
$$
\psi:\GL_m(k)\times\GL_n(k)\to G
$$
defined by
$$
\forall M\in R_1=\bigoplus_{i,j} kX_{ij}=M_{m\times n}(k)\quad
\psi{\big(}(\gamma_1,\gamma_2){\big)}(M)=\gamma_1M\gamma_2^{-1}.
$$
Then $m\neq n$ implies $G^0=G$. In case $m=n$ we have $G/G^0=\ZZ_2$ where the
other class is represented by the matrix transposition. Moreover, scheme
theoretically $G^0$ is the $k$-rational locus of the unity component of $G$.
\end{proposition}

This is Corollary 3.4 in \cite{BG1}. Here we just sketch how it follows from
Theorem \ref{t3.2}. In the case $r=1$ the determinantal ring is just the
polytopal ring corresponding to the polytope $\Delta_{m-1}\times\Delta_{n-1}$,
that is the coordinate ring of the Segre embedding
$\PP^{m-1}\times\PP^{n-1}\to\PP^{mn-1}$ and Theorem \ref{t3.2} applies. For
higher $r$ we look at the singular locus of $R$ which is exactly the coordinate
ring of the locus of matrices with rank at most $r-1$. Since the singular locus
is invariant under the automorphisms, the induction process goes through.

Now we describe those results from \cite{BG2} that are relevant in Section 4
below. At this point we need to require that the field $k$ is algebraically
closed. The case of arbitrary fields remains open.

Let $P\subset\RR^n$ be a lattice polytope of dimension $n$ and $F\subset
P$ a face. Then there is a uniquely determined retraction
$$
\pi_F:k[P]\to k[F],\ \ \pi_F(x)=0\quad \text{for}\ x\in\L_P
\setminus F.
$$
Retractions of this type will be called {\em face retractions}, and {\em facet
retractions} if $F$ is a facet. In the latter case we write $\codim(\pi_F)=1$.

Now suppose there are an affine subspace $H\subset\RR^n$
and a vector subspace $W\subset\RR^n$ with $\dim W+\dim H=n$, such that
$$
\L_P\subset\bigcup_{x\in\L_P\cap H}(x+W).
$$
(Observe that $\dim(H\cap P)=\dim H$.) The triple $(P,H,W)$ is called a {\em
lattice fibration of codimension $c=\dim W$}, whose {\em base polytope} is
$P\cap H$; its {\em fibers} are the maximal lattice subpolytopes of $(x+W)\cap
P$, $x\in\L_P\cap H$ (the fibers may have smaller dimension than $W$). $P$
itself serves as a {\em total polytope} of the fibration. If $W=\RR w$ is a
line, then we call the fibration {\em segmental} and write $(P,H,w)$ for it.
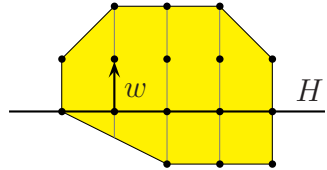
\begin{figure}[htb]
\begin{center}
\begin{pspicture}(0,0.3)(5,2.7)
\def\vertex{\pscircle[fillstyle=solid,fillcolor=black]{0.07}}
\pspolygon[fillstyle=solid,fillcolor=yellow,linewidth=0.4pt](1,1)(3,0)(5,0)(5,2)(4,3)(2,3)(1,2)
\psline[linewidth=0.4pt,linecolor=gray](2,0.5)(2,3)
\psline[linewidth=0.4pt,linecolor=gray](3,0)(3,3)
\psline[linewidth=0.4pt,linecolor=gray](4,0)(4,3)
\psline(0,1)(6,1)
\multirput(3,0)(1,0){3}{\vertex}
\multirput(1,1)(1,0){5}{\vertex}
\multirput(1,2)(1,0){5}{\vertex}
\multirput(2,3)(1,0){3}{\vertex}
\rput(5.7,1.4){$H$}
\rput(2.4,1.4){$w$}
\psline{->}(2,1)(2,1.95)
\end{pspicture}
\end{center}
\caption{A lattice segmental fibration}\label{Fig3}
\end{figure}
Note that the column structures give rise to lattice segmental fibrations in a
natural way.

For a lattice fibration $(P,H,W)$ let $L\subset\ZZ^n$ denote the subgroup
spanned by $\L_P$, and let $H_0$ be the translate of $H$ through the
origin. Then one has the direct sum decomposition
$$
L=(L\cap W)\oplus(L\cap H_0).
$$
Equivalently,
$$
\gp(S_P)=L\oplus\ZZ= \bigl(\gp(S_P)\cap W_1\bigr)\oplus
\gp(S_{P\cap H_1}).
$$
where $W_1$ is the image of $W$ under the embedding $\RR^n\to\RR^{n+1}$,
$w\mapsto (w,0)$, and $H_1$ is the vector subspace of $\RR^{n+1}$ generated by
all the vectors $(h,1)$, $h\in H$.

For a fibration $(P,H,W)$ one has the naturally associated retraction
$$
\rho_{(P,H,W)}:k[P]\to k[{P\cap H}];
$$
it maps $\L_P$ to $\L_{P\cap H}$ so that fibers are contracted to their
intersection points with the base polytope $P\cap H$.

The following is a composition of Theorems 2.2 and 8.1 in \cite{BG2}:

\begin{theorem}\label{t3.4}
\begin{itemize}
\item[(a)]
If $A$ is a retract of a polytopal algebra $k[P]$ and $\dim A\leq 2$ then $A$
is polytopal itself, i.~e. $A$ is of type either $k$ or $k[t]$ or
$k[c\Delta_1]$ for some $c\in\NN$.
\item[(b)]
If $\dim P=2$ and $f:k[P]\to k[P]$ is any idempotent endomorphism of
codimension $1$ ($f^2=f$ and $\dim\Im(f)=2$) then either
\begin{itemize}
\item[(i)] $P$ admits a lattice segmental fibration $(P,H,w)$ and $\gamma\circ
f\circ\gamma^{-1}=\iota\circ \rho_{(P,H,w)}$ for some $\gamma\in\Gamma_k(P)$
and an embedding $\iota:k[P\cap H]\to k[P]$, or
\item[(ii)] $\gamma\circ f\circ\gamma^{-1} =\iota\circ\pi_E$ for some
$\gamma\in\Gamma_k(P)$, an edge $E\subset P$ and an embedding $\iota:k[E]\to
k[P]$.
\end{itemize}
\end{itemize}
\end{theorem}

\begin{remark}\label{r3.5}
The conjectures (A) and (B) in \cite{BG2} say that both statements (a) and (b)
of Theorem \ref{t3.4} generalize to arbitrary dimensions. They should be
thought of as analogues of the fact that idempotent matrices are conjugate to
`subunit' matrices (having diagonal entries $0$ and $1$ and entries $0$
everywhere else). That one has to restrict oneself to codimension 1 idempotent
endomorphisms in the higher analogue of (b) is explained by explicit examples
in \cite[Section 5]{BG2}.

The splitting embeddings $\iota$ will be described in the next section.
\end{remark}

\section{Tame homomorphisms}

Assume we are given two lattice polytopes $P,Q\subset\RR^d$ and a homomorphism
$f:k[P]\to k[Q]$ in $\Pol(k)$. Under certain conditions there are several
standard ways to derive new homomorphisms from it.

First assume we are given a subpolytope $P'\subset P$ and a polytope
$Q'\subset\RR^n$, $d\leq n$, such that $f(k[P'])\subset k[Q']$. Then $f$ gives
rise to a homomorphism $f':k[P']\to k[Q']$ in a natural way. (Notice that we
may have $Q\subset Q'$.) Also if  $P\approx\tilde P$ and $Q\approx\tilde Q$ are
lattice polytope isomorphisms, then $f$ induces a homomorphism $\tilde
f:k[\tilde P]\to k[\tilde Q]$. We call these types of formation of new
homomorphisms {\em polytope changes}.

Now consider the situation when $\Ker(f)\cap S_P=\emptyset$. Then $f$ extends
uniquely to a homomorphism $\bar f:k[\bar S_P]\to k[\bar S_Q]$ of the
normalizations. Here $\bar S_P=\{x\in\gp(S_P)\ |\ x^m\in S_P\ \text{for some}\
m\in\NN\}$ and similarly for $\bar S_Q$. (It is well known that $k[\bar
S]=\overline{k[S]}$ for any affine semigroup $S\subset\ZZ^d$ where on the right
hand side we mean the normalization of the domain $k[S]$, \cite[Ch.\ 6]{BH}.)
This extension is given by
$$
\bar f(x)=\frac{f(y)}{f(z)},\qquad x=\frac yz,\ x\in\bar S_P,\  y\in S_P,\ z\in
S_Q.
$$
For any natural number $c$ the subalgebra of $k[\bar S_P]$ generated by the
homogeneous component of degree $c$ is naturally isomorphic to the polytopal
algebra $k[cP]$, and similarly for $k[\bar S_Q]$. Therefore, the restriction of
$\bar f$ gives rise to a homomorphism $f^{(c)}:k[cP]\to k[cQ]$. We call the
homomorphisms $f^{(c)}$ {\em homothetic blow-ups} of $f$. (Note that $k[cP]$ is
often a proper overring of the $c$th Veronese subalgebra of $k[P]$.)

One more process of deriving new homomorphisms is as follows. Assume that
homomorphisms $f,g:k[P]\to k[Q]$ are given such that
$$
\forall x\in\L_P\ \ \N(f(x))+\N(g(x))\subset Q,
$$
where $\N(-)$ denotes the Newton polytope and $+$ is the Minkowski sum in
$\RR^d$. Then we have $z^{-1}f(x)g(x)\in k[Q]$ where $z=(0,\dots,0,1)\in S_Q$.
Clearly, the assignment
$$
\forall x\in\L_P\ \ x\mapsto z^{-1}f(x)g(x)
$$
extends to a $\Pol(k)$-homomorphism $k[P]\to k[Q]$, which we denote by
$f\star g$. We call this process {\em Minkowski sum} of homomorphisms.

All the three mentioned recipes have a common feature: the new homomorphisms
are defined on polytopal algebras of dimension at most the dimension of the
sources of the old homomorphisms. As a result we are not able to really create
a non-trivial class of homomorphisms using only these three procedures. This
possibility is provided by the fourth (and last in our list) process.

Suppose $P$ is a pyramid with vertex $v$ and basis $P_0$ such that
$\L_P=\{v\}\cup\L_{P_0}$, that is $P=\join(v,P_0)$ in the terminology of
\cite{BG2, BG3}. Then $k[P]$ is a polynomial extension $k[P_0][v]$. In
particular, if $f_0:k[P_0]\to k[Q]$ is an arbitrary homomorphism and $q\in
k[Q]$ is any element, then $f_0$ extends to a homomorphism $f:k[P]\to k[Q]$
with $f(v)=q$. We call $f$ a {\em free extension} of $f_0$.

\begin{conjecture}\label{C4.1}
Any homomorphism in $\Pol(k)$ is obtained by a sequence of taking
free extensions, Minkowski sums, homothetic blow-ups, polytope changes
and compositions, starting from the identity mapping $k\to k$. Moreover,
there are normal forms of such sequences for idempotent endomorphisms.
\end{conjecture}

Observe that for general homomorphisms we do not mean that the
constructions mentioned in the conjecture are to be applied in certain
order so that we get normal forms: we may have to repeat a procedure of
the same type at different steps. However, the results mentioned in
Section 3 and Theorem \ref{t4.3} below show that for special classes of
homomorphisms such normal forms are possible.

We could call the homomorphisms obtained in the way described by Conjecture
\ref{C4.1} just {\em tame}. Then we have the {\em tame} subcategory
$\Pol(k)_{\tame}$ (with the same objects), and the conjecture asserts that
actually $\Pol(k)_{\tame}=\Pol(k)$.

\begin{remark}\label{r4.2}
(a) The correctness of Conjecture \ref{C4.1} may depend on whether or not $k$
is algebraically closed. For instance, some of the arguments in \cite{BG2} only
go through for algebraically closed fields.

(b) The current notion of tameness is weaker then the one for retractions and
surjections in \cite{BG2}. This follows from Example 5.2 and 5.3 of \cite{BG2}
in conjunction with the observation that all the explicitly constructed
retractions in \cite{BG2} are tame in the new sense.

(c) Theorems \ref{t3.2} and \ref{t3.4} can be viewed as substantial refinements
of the conjecture above for the corresponding classes of homomorphisms. Observe
that the tameness of elementary automorphisms follows from their alternative
description in Section 3. We also need the tameness of the following classes of
homomorphisms: automorphisms that map monomials to monomials, retractions of
the type $\rho_{(P,H,w)}$ and $\pi_F$ and the splitting embeddings $\iota$ as
in Theorem \ref{t3.4}. This follows from Theorem \ref{t4.3} and Corollary
\ref{c4.4} below.
\end{remark}

The next result shows that certain basic classes of morphisms in $\Pol(k)$
are tame.

\begin{theorem}\label{t4.3}
Let $k$ be a field (not necessarily algebraically closed). Then
\begin{itemize}
\item[(a)]
any homomorphism from $k[c\Delta_n]$, $c,n\in\NN$, is tame,
\item[(b)]
if $\iota:k[c\Delta_n]\to k[P]$ ($c\in\NN$) splits either $\rho_{(P,H,W)}$ for
some lattice fibration $(P,H,W)$ or $\pi_E$ for some face $E\subset P$ then
there is a normal form for representing $\iota$ in terms of certain basic tame
homomorphisms.
\end{itemize}
\end{theorem}

\begin{corollary}\label{c4.4}
For every field $k$ the homomorphisms respecting monomial structures are tame,
and we have the inclusion $\Vect_{\NN}(k)\subset\Pol_{\mathrm{tame}}(k)$, where
$\Vect_{\NN}(k)$ is the full subcategory of $\Pol(k)$ spanned by the
objects of the type $k[c\Delta_n]$, $c,n\in\NN$.
\end{corollary}

\begin{proof}[Proof of Corollary \ref{c4.4}]
The inclusion $\Vect_{\NN}(k)\subset\Pol_{\mathrm{tame}}(k)$ follows from
Theorem \ref{t4.3}(a).

Assume $f:k[P]\to k[Q]$ is a homomorphism respecting the monomial structures
and such that $\Ker(f)\cap S_P=\emptyset$. By a polytope change we can assume
$P\subset c\Delta_n$ for a sufficiently big natural number $c$, where $n=\dim
P$ and $\Delta_n$ is taken in the lattice $\ZZ L_P$. In this situation there is
a bigger lattice polytope $Q'\supset Q$ and a unique homomorphism
$g:k[c\Delta_n]\to k[Q']$ for which $g|_{\L_P}= f|_{\L_P}$. By Theorem
\ref{t4.3}(a) $f$ is tame.

Consider the situation when the ideal $I=(\Ker(f)\cap S_P)k[P]$ is a nonzero
prime monomial ideal and there is a face $P_0\subset P$ such that
$\Ker(f)\cap\L_{P_0}=\emptyset$ and $f$ factors through the face projection
$\pi:k[P]\to k[P_0]$, that is $\pi(x)=x$ for $x\in\L_{P_0}$ and $\pi(x)=0$ for
$x\in\L_P\setminus\L_{P_0}$. In view of the previous case we are done once the
tameness of face projections has been established.

Any face projection is a composite of facet projections. Therefore we can
assume that $P_0$ is a facet of $P$. Let $(\RR P)_+\subset\RR P$ denote the
halfspace that is bounded by the affine hull of $P_0$ and contains $P$. There
exists a unimodular (with respect to $\ZZ\L_P$) lattice simplex
$\Delta\subset(\RR P)_+$ such that $\dim\Delta=\dim P$, the affine hull of
$P_0$ intersects $\Delta$ in one of its facets and $P\subset c\Delta$ for some
$c\in\NN$. But then $\pi$ is a restriction of the corresponding facet
projection of $k[c\Delta]$, the latter being a homothetic blow-up of the
corresponding facet projection of the polynomial ring $k[\Delta_n]$ --
obviously a tame homomorphism.
\end{proof}

\begin{proof}[Proof of Theorem \ref{t4.3}(a)]
We will use the notation
$\{x_0,\ldots,x_n\}=\L_{\Delta_n}$. Any lattice point $x\in c\Delta_n$ has a
unique representation $x=a_0x_0+\cdots+a_nx_n$ where the $a_i$ are nonnegative
integer numbers satisfying the condition $a_0+\cdots+a_n=c$. The numbers $a_i$
are the {\it barycentric coordinates} of $x$ in the $x_i$.

Let $f:k[c\Delta_n]\to k[P]$ be any homomorphism.

First consider the case when one of the points from $\L_{c\Delta_n}$ is
mapped to $0$. In this situation $f$ is a composite of facet projections
and a homomorphism from $k[c\Delta_m]$ with $m<n$.
As observed in the proof of Corollary \ref{c4.4} facet projections are
tame. Therefore we can assume that none of the $x_i$ is mapped to $0$.
By a polytope change we can also assume
$\L_P\subset\{X_1^{a_1}\cdots X_r^{a_r}Y^bZ\ |\ a_i,b\geq0\}$, $r=\dim P-1$.

Consider the polynomials $\phi_x=Z^{-1}f(x)\in k[X_1,\ldots,X_r,Y]$,
$x\in\L_{c\Delta_n}$. Then the $\phi_x$ are subject to the same binomial
relations as the $x$. One the other hand the multiplicative semigroup
$k[X_1,\ldots,X_r,Y]\setminus\{0\}/k^*$ is a free commutative semigroup
and, as such, is an inductive limit of free commutative semigroups
of finite rank. Therefore, by Lemma \ref{l4.4} below
there exist polynomials $\psi,\eta_i\in k[X_1,\ldots,X_r,Y]$, $i\in[0,n]$,
and scalars $t_x\in k^*$, $x\in\L_{c\Delta_n}$,
such that $\phi_x=t_x\psi\eta_0^{a_0}\cdots\eta_n^{a_n}$
where the $a_i$ are the barycentric coordinates of $x$. Clearly, $t_x$ are
subject to the same binomial relations as the $x\in\L_{c\Delta_n}$. Therefore,
after the normalizations $\eta_i\mapsto t_{x_i}\eta_i$ ($i\in[0,n]$) we get
$\phi_x=\psi\eta_0^{a_0}\cdots\eta_n^{a_n}.$ But the latter equality can
be read as follows: $f$ is obtained by a polytope change applied to
$\Psi\star\Theta^{(c)}$, where
\begin{itemize}
\item[(i)]
$\Psi:k[c\Delta_n]\to k[Q]$,\ $\Psi(x)=\psi Z$,\ $x\in\L_{c\Delta_n}$,
\item[(ii)]
$\Theta:k[\Delta_n]\to k[Q]$,\ $\Theta(x_i)=\eta_iZ$,\ $i\in[0,n]$,
\end{itemize}
and $Q$ is a sufficiently large lattice polytope so that it contains all
the relevant lattice polytopes. Now $\Psi$ is tame because it can be
represented as the composite map
$$
k[c\Delta_n]\xrightarrow{\L_{c\Delta_n}\to t}k[t]\xrightarrow
{t\mapsto\psi Z}k[Q]
$$
(the first map is the $c$th homothetic blow-up of $k[\Delta_n]\to k[t]$,
$x_i\mapsto t$ for all $i\in[0,c]$) and $\Theta$ is just a free extension
of the identity embedding $k\to k[Q]$.

(b) First consider the case of lattice segmental fibrations.

Consider the rectangular prism $\Pi=(c\Delta_n)\times(m\Delta_1)$.
By a polytope change (assuming $m$ is sufficiently large) we can assume that
$P\subset\Pi$ so that $H$ is parallel to $c\Delta_n$:
\begin{figure}[hbt]
\begin{pspicture}(-0.5,-0.5)(3.5,4.5)
\def\vertex{\pscircle[fillstyle=solid,fillcolor=black]{0.07}}
\pspolygon[fillstyle=solid,fillcolor=yellow,linewidth=0.7pt,linecolor=black](0,1)(1,0)%
   (2,0)(3,1)(3,3)(1,4)(0,2)
\psset{linewidth=0.4pt,linecolor=gray}
\multirput(0,0)(0,1){5}{\psline(0,0)(3,0)}
\multirput(0,0)(1,0){4}{\psline(0,0)(0,4)}
\psset{linewidth=1pt,linecolor=black}
\multirput(0,0)(0,1){5}{\vertex}
\multirput(1,0)(0,1){5}{\vertex}
\multirput(2,0)(0,1){5}{\vertex}
\multirput(3,0)(0,1){5}{\vertex}
\psline(0,1)(0,2)
\psline(1,0)(2,0)
\psline(3,1)(3,3)
\psline(-1,1)(4,1)
\rput[bl](-1.5,1.8){$m\Delta_1$}
\rput[bl](1.1,-0.8){$c\Delta_n$}
\rput(1.5,2.5){$\Pi$}
\rput[bl](3.4,1.2){$H$}
\end{pspicture}
\caption{}\label{fig4}
\end{figure}
The lattice point $(x,b)\in\Pi$ will be identified with the monomial
$X_1^{a_1}\cdots X_n^{a_n}Y^bZ$ whenever we view it as a monomial in $k[\Pi]$,
where the $a_i$ are the corresponding barycentric coordinates of $x$ (see the
proof of (a) above). (In other words, the monomial $X_1^{a_1}\cdots X_n^{a_n}$
is identified with the point $x=(c-a_1-\cdots-a_n)x_0+a_1x_1+\cdots+a_nx_n\in
c\Delta_n$.)

Assume $A:k[c\Delta_n]\to k[m'\Delta_1]$ is a homomorphism of the type
$A(x_i)=a Z$, $i\in[0,c]$ for some $a\in k[Y]$ satisfying the
condition $a(1)=1$. Consider any homomorphism $B:k[\Delta_n]\to k[\Pi']$,
$\Pi'=\Delta_n\times(m'\Delta_1)$ that splits the projection
$\rho':k[\Pi']\to k[\Delta_n]$,
$\rho'(ZY^b)=Z$ and $\rho'(X_iY^bZ)=X_iZ$ for $i\in[1,n]$, $b\in[0,m']$.
The description of such  homomorphisms is clear -- they are exactly the
homomorphisms $B$ for which
\begin{itemize}
\item[]
$B(x_0)=B_0\in Z+(Y-1)(Zk[Y]+X_1Zk[Y]+\cdots+X_nZk[Y])$,\smallskip
\item[]
$B(x_i)=B_i\in X_iZ+(Y-1)(Zk[Y]+X_1Zk[Y]+\cdots+X_nZk[Y])$,
$i\in[1,n]$,\smallskip
\item[]
$\deg_YB_i\leq m'$ for all $i\in[0,n]$.
\end{itemize}
Clearly, all such $B$ are tame.

In case $m\geq\max\{m'+c\deg_YB_i\}_{i=0}^n$ we have the homomorphism
$A\star B^{(c)}:k[c\Delta_n]\to k[\Pi]$ which
obviously splits the projection $\rho:k[\Pi]\to
k[c\Delta_n]$ defined by
$$
\rho(X_1^{a_1}\cdots X_n^{a_n}Y^bZ)=X_1^{a_1}\cdots X_n^{a_n}Z.
$$
Assume $\iota$ splits $\rho_{(P,H,w)}$. Since the Newton
polytope of a product is a Minkowski sum of the Newton
polytopes of the factors, we get: the polynomials $\psi$ and
$\eta_i$, mentioned in the proof of (a), that correspond to
$\iota$, satisfy the conditions: $\psi\in k[Y]$ and
$\eta_i\in k[Y]+X_1k[Y]+\cdots+X_nk[Y]$.
Is is also clear that upon evaluation at $Y=1$ we get
$\psi(1),\eta_i(X_1,\ldots,X_n,1)\in k^*$, $i\in[0,n]$.
Therefore, after the normalizations $\psi\mapsto\psi^{-1}(1)\psi$,
$\eta_i\mapsto\eta_i(X_1,\ldots X_n,1)^{-1}\eta_i$
we conclude that $\iota$ is obtained by a polytope change applied to
$A\star B^{(c)}$ as above (with respect to $a=\psi$, $B_0=\eta_iZ$,
$i\in[0,n]$.

For a lattice fibration $(P,H,W)$ of higher codimension similar arguments show
that $\iota$ is obtained by a polytope change applied to $A\star B^{(c)}$,
where
\begin{itemize}
\item[]
$B$ is a splitting of a projection of the type $\rho_{(P',H',W')}$ such
that the base polytope $P'\cap H'$ is a unit simplex and\smallskip
\item[]
$A$ is a homomorphism defined by a single polynomial whose Newton
polytope is parallel to $W'$.
\end{itemize}

We skip the details for splittings of face projections and only remark that
similar arguments based on Newton polytopes imply the following. All such
splittings are obtained by polytope changes applied to $A\star B^{(c)}$
where $B$ is a splitting of a face projection onto a polynomial ring and
$A$ is again defined by a single polynomial.
\end{proof}

\begin{lemma}\label{l4.4}
Assume we are given an integral affine mapping
$\alpha:c\Delta_n\to\RR_+^d$ for some natural numbers $c$, $n$ and $d$.
Then there exists an element $v\in\ZZ_+^d$ and a integral affine
mapping $\beta:\Delta_n\to\RR^d_+$ such that $\alpha=v+c\beta$.
\end{lemma}
\begin{proof}
Assume $\alpha(cx_i)=(a_{i1},\ldots,a_{id})$, $i\in[0,n]$
(the $x_i$ as above). Consider the vector
$$
v={\big(}\min\{a_{i1}\}_{i=0}^n,\ldots,\min\{a_{id}\}_{i=0}^n{\big)}.
$$
It suffices to show that all the vectors $\alpha(cx_i)-v$ are $c$th multiples
of integral vectors. But for any index $l\in[1,d]$ the $l$th component of
either $\alpha(cx_i)-v$ or $\alpha(cx_j)-v$ for some $j\not= i$ is zero. In the
first case there is nothing to prove and in the second case the desired
divisibility follows from the fact that
$\alpha(cx_i)-\alpha(cx_j)=(\alpha(cx_i)-v)-(\alpha(cx_j)-v)$ is a $c$th
multiple of an integral vector (because $\alpha$ is integral affine).
\end{proof}

\begin{remark}
Theorems \ref{t3.2}, \ref{t3.4} and \ref{t4.3} provide a possibility for
computing the $\Gamma_k(P)$-variety of idempotent endomorphisms
$\Idemp(k[P])$  for a polygon $P$ and an algebraically closed field $k$
($\Gamma_k(P)$ acts by conjugation): the orbits of
codimension 1 idempotent endomorphisms are naturally associated to the
segmental fibration structures and edges of $P$, and all codimension 2
idempotent endomorphisms factor through $k[t]$.
\end{remark}

\end{document}